\documentclass[12pt,twoside,a4paper]{article}
\usepackage{amsmath,amssymb,amsthm}

\usepackage{graphicx,psfrag}

\newtheorem{Thm}{Theorem}[section]
\newtheorem{Lem}[Thm]{Lemma}
\newtheorem{Pro}[Thm]{Proposition}

\theoremstyle{definition}

\theoremstyle{remark}
\newtheorem{Rem}[Thm]{Remark}

\newcommand{\R}{\mathbb{R}}
\newcommand{\Z}{\mathbb{Z}}
\newcommand{\N}{\mathbb{N}}

\newcommand{\cT}{\mathcal{T}}

\newcommand{\al}{\alpha}

\newcommand{\ga}{\gamma}

\newcommand{\de}{\delta}

\newcommand{\ep}{\varepsilon}

\newcommand{\si}{\sigma}

\newcommand{\la}{\lambda}
\newcommand{\La}{\Lambda}
\renewcommand{\phi}{\varphi}

\newcommand{\dist}{\operatorname{dist}}
\newcommand{\diam}{\operatorname{diam}}

\newcommand{\hyp}{\operatorname{H}}
\newcommand{\id}{\operatorname{id}}

\newcommand{\pr}{\operatorname{pr}}

\newcommand{\const}{\operatorname{const}}
\newcommand{\es}{\emptyset}
\renewcommand{\d}{\partial}
\newcommand{\di}{\d_{\infty}}
\newcommand{\set}[2]{\{#1:\,\text{#2}\}}
\newcommand{\sm}{\setminus}
\newcommand{\sub}{\subset}

\newcommand{\ov}{\overline}

\newcommand{\md}{\!\mod}

\begin{document}

\title{Embedding of hyperbolic spaces in the product of trees}
\author{Sergei Buyalo\footnote{Supported by RFFI Grant
02-01-00090, CRDF Grant RM1-2381-ST-02 and
SNF Grant 20-668 33.01}
\ \& Viktor Schroeder\footnote{Supported by Swiss National Science 
Foundation}}

\date{}
\maketitle

\begin{abstract} We show that for each
$n\ge 2$
there is a quasi-isometric embedding of the hyperbolic
space
$\hyp^n$
in the product
$T^n=T\times\dots\times T$
of
$n$
copies of a (simplicial) metric tree
$T$.
On the other hand, we prove that there is no quasi-isometric
embedding
$\hyp^2\to T\times\R^m$
for any metric tree
$T$
and any
$m\ge 0$.
\end{abstract}

\section{Introduction}
 Recall that a map
$f:X\to Y$
between metric spaces is called a large scale uniform embedding if
$$\phi_1(|x-x'|)\le|f(x)-f(x')|\le\phi_2(|x-x'|)$$
for some functions
$\phi_1$, $\phi_2:[0,\infty)\to[0,\infty)$
tending to infinity and all
$x$, $x'\in X$. The map
$f$
is called quasi-isometric, if one can take linear functions as
$\phi_1$, $\phi_2$,
$\phi_i(t)=l_it+m_i$, $i=1,2$.

We denote by
$\hyp^n$
the real hyperbolic space of dimension
$n$
and of curvature
$-1$.

\begin{Thm}\label{Thm:main1} For each 
$n\ge 2$
there is a quasi-isometric embedding
$$f:\hyp^n\to T^n=T\times\dots\times T,$$
where
$T$
is a homogeneous simplicial metric tree, whose edges all
have length 1.
\end{Thm}

\begin{Rem}\label{Rem:infval} Every vertex of the tree 
$T$
from Theorem~\ref{Thm:main1} has infinite (countable)
valence, i.e., it is adjacent to infinitely many edges.
In particular, 
$T$
is not locally compact. For
$n=2$
there is a better result \cite{DS}, saying that the hyperbolic
plane
$\hyp^2$
can be quasi-isometrically embedded in the product of
two locally compact metric trees.
\end{Rem}

\begin{Rem}\label{Rem:asdim} There is a general embedding
result \cite{Dr} according to which every metric space of
bounded geometry, whose asymptotic dimension
$\le n$,
admits a large scale uniform embedding into the product
of
$n+1$
locally compact metric trees.
The hyperbolic space
$\hyp^n$
has bounded geometry and its asymptotic dimension equals
$n$.
Thus our Theorem~\ref{Thm:main1} is stronger than the
Dranishnikov's result applied to
$\hyp^n$
w.r.t. the number of trees needed for an embedding and
the quality of embeddings: we construct quasi-isometric 
embeddings. On the other hand, it is weaker w.r.t. finiteness 
properties of the target trees.
\end{Rem}

One can ask whether it is possible to embed 
$\hyp^n$
quasi-isometrically in the product of less than
$n$
metric trees. To make this question nontrivial, one should
stabilize the product by an additional factor which has
arbitrarily large dimension and small growth rate, e.g.,
by
$\R^m$.
It easily follows from results of our previous paper
\cite{BS1} that there is no quasi-isometric
embedding
$\hyp^n\to X$,
$X=T_1\times\dots\times T_p\times\R^m$,
for any
$p\le n-2$ 
and
$m\ge 0$.
For the projection
$X\to T_1\times\dots\times T_p$
defines a subexponential foliation of
$X$ 
of rank
$p=\dim (T_1\times\dots\times T_p)$,
therefore, the subexponential corank of
$X$
is
$\le p$,
and by the main result of \cite{BS1}, the existence of
$\hyp^n\to X$
implies
$n-1\le p$.
In fact, Theorem~\ref{Thm:main1}
is optimal w.r.t. the number of trees in the product stabilized
by
$\R^m$.
We prove here that this is true for
$n=2$.
The general case is considered in a fortcoming paper \cite{BS2}.

\begin{Thm}\label{Thm:main2} There is no quasi-isometric
embedding
$\hyp^2\to T\times\R^m$
for any metric tree
$T$
and
$m\ge 0$.
\end{Thm}

\medskip
{\em Acknowledgment.} The first author is happy to express his
deep gratitude to the University of Z\"urich for the support,
hospitality and excellent working conditions while writing the 
paper.

\section{Proof of Theorem~\ref{Thm:main1}}

\subsection{Idea of the embedding}

We describe the idea of the embedding for the case
$n=2$.
We write
$\hyp^2$
in horospherical coordinates
$\hyp^2=\R\times\R$
such that the sets
$\{t\}\times\R$
are horocycles. Consider the integer horocycle
$h_i=\{i\}\times\R$
with intrinsic metric isometric to the real line. The
canonical projection
$\pi:h_i\to h_{i-1}$
is a homothety. We choose an integer
$p\ge 5$ 
and assume (after scaling the metric of
$\hyp^2$
suitable) that homothery factor of
$\pi$
is
$1/p$.

Consider on each horocycle
$h_i$
in a periodic way intervals
$Q_{ij}$, $j\in\Z$,
all of length
$(1-\frac{2}{p})<1$
such that the gap between two neighboring intervals is
$2/p<1/2$.
\begin{figure}[htbp]
\centering
\includegraphics[width=1.0\columnwidth]{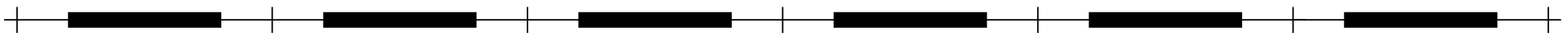}
\caption{$p=5$}\label{Fi:2}
\end{figure}

It is not difficult to arrange these intervals in a way that
\begin{itemize}
\item[(i)] the projection
$\pi(Q_{ij})\sub h_{i-1}$
on an interval
$Q_{ij}$
is either contained completely in some interval
$Q_{i-1,j'}$
or completely in a gap;
\item[(ii)] for every
$Q_{ij}$
there exists a
$k>0$
such that
$\pi^k(Q_{ij})$
is contained in some interval
$Q_{i-k,j'}\sub h_{i-k}$.
\end{itemize}

The intervals now define a tree
$T$:
the vertices are the intervals
$Q_{ij}$.
The vertex
$Q_{ij}$
is connected by an edge with
$Q_{i-k,j'}$,
where
$k=k(i,j)$
and
$Q_{i-k,j'}$
is the smallest integer and the interval according to (ii).
The map
$f:\hyp^2\to T$
is defined by associating to a point
$x$
an interval
$Q_{ij}$
with minimal distance to
$x$; $f$
is Lipschitz on a large scale due to (i).

Next we define a second tree
$T'$
in the same way using now intervals
$Q_{ij}'$
such that for every
$i$, $\cup_{j\in\Z}(Q_{ij}\cup Q_{ij}')=h_i$,
i.e., the intervals
$Q_{ij}'$, $j\in\Z$,
cover the gap of the intervals
$Q_{ij}$.
Finally, we will show that
$(f,f'):\hyp^2\to T\times T'$
is quasi-isometric.

\medskip
For convenience of notations, we shift the dimension by 1,
and construct a quasi-isometric embedding
$\hyp^{n+1}\to T^{n+1}$
assuming that
$n\ge 1$.

\subsection{Construction of the target tree
$T$}\label{subsect:targetree}

To construct the target tree
$T$,
we consider the unit cube
$I^n\sub\R^n$,
where
$I=[0,1]\sub\R$.
We fix
$a\in(0,1)$
such that
$p:=\frac{2}{1-a}$
is an integer,
$p\in\N$,
and moreover
$\frac{1}{p-1}+\frac{1}{p}<\frac{1}{n+1}$.
In particular,
$p>2(n+1)$.

Consider the subsegment
$J\sub I$
of the length
$a$
centered at the middle of
$I$,
i.e., at
$1/2$.
Now, the middle subcube
$A=J^n\sub I^n$
will play the role of a template for the vertices of
$T$.
The condition
$\frac{1}{p-1}+\frac{1}{p}<\frac{1}{n+1}$
will be used while constructing 
$n+1$
appropriate copies of 
$T$,
see sect.~\ref{subsect:color}.

\subsubsection{Definition of the vertices of
$T$}\label{subsubsect:defvert}

Using the action of the integer lattice
$\Z^n\sub\R^n$
by shifts on
$\R^n$
we define the set
$$Q_0=\bigcup_{\ga\in\Z^n}\ga A.$$
Note that
$Q_0$
is a disconnected subset of
$\R^n$, 
every connected component
$\ga A$, $\ga\in\Z^n$,
of which is a cube of diameter
$a\sqrt{n}$.
We call the connected components of the set
$Q_0$
the vertices of
$Q_0$,
and we shall identify them with the vertices 
of the level 0 of the tree
$T$.

To define the vertices of the level 1, we apply the
following procedure. We subdivide the segment
$I$
into
$p$
equal subsegments of the length
$1/p$, 
so that
$p-2$
of them cover the middle subsegment
$J$,
and the remaining two cover its complement in
$I$.
This subdivision induces the subdivision of the cube
$I^n$
into
$p^n$
congruent and parallel subcubes. There is a natural labeling
of these subcubes 
$I_l^n$, $l\in L$,
by the set
$L:=\{1,\dots,p\}^n$,
and for every
$l\in L$
the canonical homothety
$h_l:I^n\to I_l^n$
with the coefficient
$\la=1/p$
maps the middle subcube
$A\sub I^n$
onto the subcube
$A_l=h_l(A)\sub I_l^n$.
Now, we define
$$Q_1=\bigcup_{\ga\in\Z^n}\bigcup_{l\in L}\ga A_l.$$
Note that
$Q_1$
is a disconnected subset of
$\R^n$, 
every connected component
$\ga A_l$, $\ga\in\Z^n$, $l\in L$,
of which is a cube of diameter
$\la a\sqrt{n}$.
We call the connected components of
$Q_1$
the vertices of
$Q_1$. 

\begin{figure}[htbp]
\centering
\includegraphics[width=1.0\columnwidth]{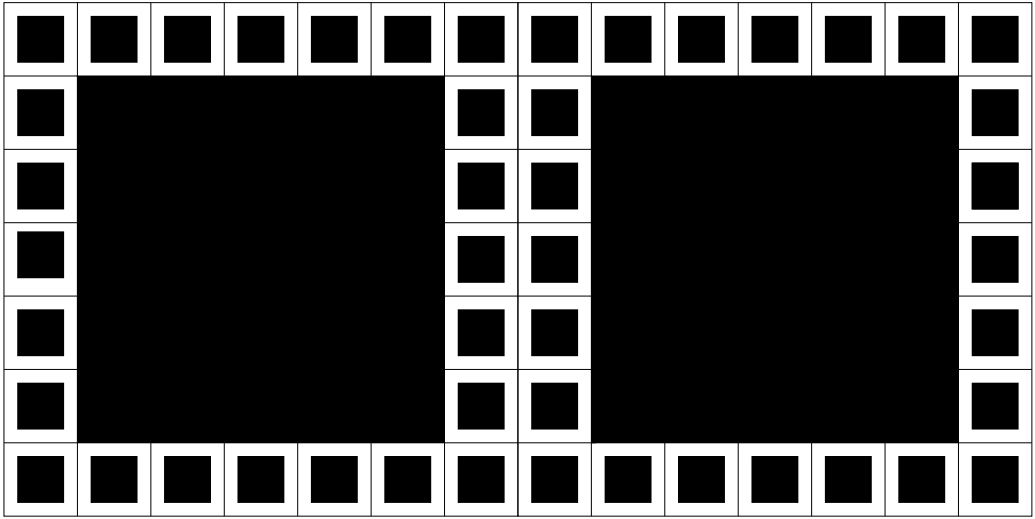}
\caption{some small cubes of the next level are hidden behind
the large black cubes; here $p=7$}\label{Fi:1}
\end{figure}

Any vertex of
$Q_1$
is either separated from every vertex of
$Q_0$
by the distance at least
$\la^2$,
or it lies inside of some vertex of
$Q_0$
being separated from its boundary by the distance at least
$\la^2$.
This property is of key importance in what follows and it
is called {\em the separation property of}
$Q_0\cup Q_1$.
We shall identify the vertices of
$Q_1$
with the vertex set of level 1
of the tree
$T$.

To define the vertex set of any level 
$k\ge 1$, 
we apply repeatedly the described procedure. Namely,
consider the set
$L$
as an alphabet, and let
$W_k$
be the set of words of length
$k$
in the alphabet
$L$,
i.e., each
$w\in W_k$
is a sequence of
$k$
letters from
$L$.
In particular,
$W_0=\es$
and
$W_1=L$.
For each
$w\in W_k$, $w=l_1\dots l_k$,
we define the homothety
$h_w:I^n\to I^n$
as the composition
$h_w=h_{l_1}\circ\dots\circ h_{l_k}$, $h_{\es}=\id$.
This is a homothety with the coefficient
$\la^k$.
We let
$A_w=h_w(A)$
be a subcube in
$I^n$,
and note that
$A_w\sub A_{w'}$,
where
$w'\in W_{k-1}$,
if and only if 
$w'\sub w$
is the initial subword and
every coordinate of the last letter
$l_k\in\{1,\dots,p\}^n$
is different from
$1$
and
$p$.
Now, we define
$$Q_k=\bigcup_{\ga\in\Z^n}\bigcup_{w\in W_k}\ga A_w.$$
Again, for each
$k\ge 1$,
$Q_k$
is a disconnected set of
$\R^n$, 
every connected component
$\ga A_w$, $\ga\in\Z^n$, $w\in W_k$,
of which is a cube of diameter
$\la^k a\sqrt{n}$.
We call the connected components of
$Q_k$
the vertices of
$Q_k$. 

The separation property of the set
$Q^+=\cup_{k\ge 0}Q_k$
is the following. For each
$0\le k'<k$
any vertex of
$Q_k$
either is separated from every vertex of
$Q_{k'}$
by the distance at least
$\la^{k+1}$
or it lies inside of some vertex of
$Q_{k'}$
being separated from its boundary by the distance at least
$\la^{k+1}$.
This immediately follows from self-similarity of our 
construction. More precisely, for some vertices
$\ga A_w\sub Q_k$, $\ga'A_{w'}\sub Q_{k'}$, 
we have
$\ga A_w\sub\ga'A_{w'}$
if and only if
$\ga=\ga'$
and
$w'$
is an initial subword of
$w$,
and the first letter from
$w\sm w'$, $l_{k+1}$,
has all coordinates different from
$1$
and
$p$.
The vertices of
$Q^+$
we shall identify with the vertices of all levels 
$\ge 0$
of the tree
$T$.

To define the vertex set
$Q_{-k}$
for
$k>0$
we take the vector
$\eta=\frac{1}{p}\theta\in\R^n$,
where
$\theta:=\{1,\dots,1\}\in\R^n$,
and consider the homothety
$H:\R^n\to\R^n$,
$$H(x)=p(x-\eta),\quad x\in\R^n.$$
Then we put
$Q_{-k}:=H^k(Q_0)$
for each
$k>0$.
Again, 
$Q_{-k}$
is a disconnected subset of
$\R^n$, 
every connected component of which is a cube of diameter
$\la^{-k}a\sqrt{n}$.

Note that
$\eta_0=\frac{1}{p-1}\theta$
is the unique fixed point for
$H$, $H(\eta_0)=\eta_0$.
Furthermore, since
$p\ge 3$,
we have
$\frac{1}{p}<\frac{1}{p-1}<1-\frac{1}{p}$,
hence,
$\eta_0$
is an interior point of the cube
$A$.
Then for
$\de:=\dist(\eta_0,\d A)>0$
we have
$$\dist(\eta_0,\d H^k(A))=p^k\de\to\infty$$
as
$k\to\infty$.

Therefore, given a vertex 
$v$
of
$Q_{k'}$, $k'\in\Z$,
the vertex 
$H^k(A)$
of
$Q_{-k}$
cover  
$v$
for all sufficiently large
$k$.
This property will provide connectedness of the tree
$T$.
Furthermore, the set
$Q=\cup_{k\in\Z}Q_k$
has the separation property exactly as it is stated
above for any
$k$, $k'\in\Z$
with
$k'<k$
(for more detail see the proof of Proposition~\ref{Pro:covertree}
below).

\subsubsection{Definition of the edges of
$T$}

Let
$V_k$
be the set of the vertices of
$Q_k$, $k\in\Z$.
We define the vertex set
$V$
of the tree
$T$
as the union
$V=\cup_{k\in\Z}V_k$.
Two vertices
$v\in V_k$, $v'\in V_{k'}$
are connected by an edge in
$T$
if and only if
$k\neq k'$,
say
$k'<k$,
$v\sub v'$
(considered as cubes in
$\R^n$),
and 
$k'$
is minimal with this property. This defines a graph
$T$
with the vertex set
$V$
and the edge set
$E$.

\begin{Lem}\label{Lem:tree} The graph
$T$
is a tree, i.e., it is connected and has no
circuit.
\end{Lem}

\begin{proof} 
$T$
has no circuit because every vertex
$v\in V_k$
is connected with at most one vertex from
$V_{k'}$
for every
$k$, $k'\in\Z$, $k'<k$,
since the vertices of
$Q_{k'}$
are separated subsets of
$\R^n$.

By the property of the homothety 
$H:\R^n\to\R^n$
indicated above, given two vertices
$v'\in V_{k'}$, $v''\in V_{k''}$
the vertex 
$v=H^k(A)$
of
$V_{-k}$
cover  
$v'$, $v''$
for all sufficiently large
$k$;
$v'$, $v''\sub v$.
By the separation property of
$Q$,
if
$v'\sub u$
for some vertex
$u\in V_m$, $m>k$,
then
$m<k'$
and
$u\sub v$.
Therefore, the vertices
$v'$, $v''$
are connected by paths in
$T$
with the vertex
$v$.
Hence,
$T$
is connected.
\end{proof}

\subsection{Construction of colored copies of
$T$}\label{subsect:color}
Here we describe 
$n+1$
colored copies of
$T$, $T_c$, $c\in C$,
for which the target space
$X$
for the embedding
$\hyp^{n+1}\to X$
will be
$X=\prod_{c\in C}T_c$.
As the set of the colors we take the cyclic group
$C=\Z/(n+1)\Z$
of order
$n+1$,
and color the tree
$T$
by its unit,
$T=T_0$.
For each
$c\in C$
the vertex set
$V_c$
of
$T_c$
will be the union
$V_c=\cup_{k\in\Z}V_{c,k}$,
where every
$V_{c,k}$
is the set of cubes in
$\R^n$.

\begin{Pro}\label{Pro:covertree} For each
$c\in C$
there is a copy
$T_c$
of the tree
$T$
with the vertex set
$V_c=\cup_{k\in\Z}V_{c,k}$
such that the set
$V_c$
of cubes in
$\R^n$
satisfies the separation property, and for every
$k\in\Z$
the union
$\cup_{c\in C}V_{c,k}$
covers
$\R^n$.
\end{Pro}

The proof is based on Lemma~\ref{Lem:covershift} below.
Identifying the opposite
$(n-1)$-faces
of the cube
$I^n$,
we obtain
$n$-torus
$P^n=I^n/\sim$,
and we consider the set
$A=J^n\sub I^n$
as a subset of
$P^n$.
Recall that
$p>2(n+1)$,
see sect.~\ref{subsect:targetree}.

\begin{Lem}\label{Lem:covershift} Consider the diagonal
action of the group
$C$
on
$P^n$,
$(c,x)\mapsto x+c\nu\md\Z^n$, 
where
$\nu=\frac{1}{n+1}\theta$.
Then the shifts of
$A$
cover
$P^n$,
$\cup_{c\in C}c(A)=P^n$.
\end{Lem}

\begin{proof} We put
$Z=P^n\sm A$.
Since
$c(Z)=P^n\sm c(A)$
for every
$c\in C$,
it suffices to show that
$$G:=\bigcap_{c\in C}c(Z)=\es.$$
Assume that it is not the case.Then
$G$
contains the
$C$-orbit
of some
$x\in Z$,
since
$c(G)=G$.
The point
$x=(x_1,\dots,x_n)$
has a coordinate which lies in
$$I_0=[0,1/p)\cup(1-1/p,1].$$
Without loss of generality, we may assume that
$x_1\in I_0$.
Since
$c(x)\in Z$, $x$
has a coordinate lying in
$$I_n=\left(\frac{n}{n+1}-\frac{1}{p},
  \frac{n}{n+1}+\frac{1}{p}\right),$$
and we may assume that
$x_n\in I_n$,
because
$I_0\cap I_n=\es$
due to the condition
$p>2(n+1)$. 
Similarly, since
$c^i(x)\in Z$,
we find that
$$x_{n-i}\in I_{n-i}=\left(\frac{n-i}{n+1}-\frac{1}{p},
  \frac{n-i}{n+1}+\frac{1}{p}\right)$$
for all
$i=0,1,\dots,n-1$,
in particular,
$x_1\in I_1$.
This is a contradiction, because
$I_0\cap I_1=\es$.
\end{proof}

We fix a universal covering
$\pi:\R^n\to P^n$,
and for
$c\in C$
we put
$Q_{c,0}=\pi^{-1}(c(A))$.
Clearly,
$Q_{c,0}=Q_0+c\nu$,
where
$\nu=\frac{1}{n+1}\theta\in\R^n$
and
$C$
is identified with the set
$\{0,1,\dots,n\}$.
The cubes of
$Q_{c,0}$
form the vertex set
$V_{c,0}$
of the
$0$-level
of
$T_c$.
It follows from Lemma~\ref{Lem:covershift} that
$\cup_{c\in C}V_{c,0}=\R^n$.

Given a letter
$l\in L=\{1,\dots,p\}^n$,
we assume that the homothety
$h_l:I^n\to I_l^n$
(see sect.~\ref{subsubsect:defvert}) is canonically extended to the homothety
$\R^n\to\R^n$
for which we use the same notation.

\begin{Lem}\label{Lem:shiftalph} For every letters
$l$, $l'\in L$
and each color
$c\in C$,
the sets
$h_l(Q_{c,0})$
and
$h_{l'}(Q_{c,0})$
coincide,
$h_l(Q_{c,0})=h_{l'}(Q_{c,0})$.
\end{Lem}

\begin{proof} The cube
$c(A)$
is obtained from
$A$
by a shift,
$c(A)=A+c\nu$.
Thus the cubes
$h_l(c(A))$
and
$h_{l'}(c(A))$
are obtained from
$A_l=h_l(A)$
and
$A_{l'}=h_{l'}(A)$
respectively
by one and the same shift of
$R^n$.
The homothety
$h_{l'}$
can be obtained by composing the homothety
$h_l$
with the shift of
$\R^n$,
which moves the cube
$I_l^n$
to the cube
$I_{l'}^n$.
Since the shifts of
$\R^n$
commute, the last one moves
$h_l(c(A))$
to
$h_{l'}(c(A))$.
This easily implies the claim.
\end{proof}

Given
$k\ge 0$, $c\in C$,
we define
$Q_{c,k}:=h_w(Q_{c,0})$
for some word 
$w\in W_k$.
By Lemma~\ref{Lem:shiftalph}, this is independent
of
$w$.
The cubes of
$Q_{c,k}$
form the vertex set
$V_{c,k}$
of the
$k$-level
of
$T_c$.
It follows that
$\cup_{c\in C}V_{c,k}=\R^n$
for each
$k\ge 0$.

Furthermore, for each color
$c\in C$,
the set
$Q_c^+=\cup_{k\ge 0}Q_{c,k}$
has the separation property because every
$Q_{c,k}$, $k\ge 0$
can be obtained from the cube
$c(A)$
by the self-similarity maps
$\set{h_w}{$w\in W_k$}$
and then applying the action of
$\Z^n$.

Given a color
$c\in C$,
we define the set
$Q_{c,k}$
for negative levels as
$$Q_{c,k}:=H^{-k}(Q_{c,0}),\quad k<0,$$
where the homothety
$H:\R^n\to\R^n$
is defined at the end of sect.~\ref{subsubsect:defvert},
$H(x)=p(x-\eta)$.
The cubes of
$Q_{c,k}$
form the vertex set
$V_{c,k}$
of the
$k$-level
of 
$T_c$.
It follows that
$\cup_{c\in C}V_{c,k}=\R^n$
for each
$k<0$.

Now, we define the vertex set of
$T_c$
as
$V_c=\cup_{k\in\Z}V_{c,k}$.
The edges of
$T_c$
are defined by the same condition as for the tree
$T=T_0$.
To prove that
$T_c$
is connected, we need

\begin{Lem}\label{Lem:colorcov} The fixed point of
$H$,
$\eta_0=\frac{1}{p-1}\theta$,
lies in the interior of the cube
$c(A)$
(taken
$\md\Z^n$)
for every
$c\in C$.
\end{Lem}

\begin{proof} The point
$\eta_0(c)=\eta_0+(n+1-c)\nu$
lies in the interior of the cube
$A\md\Z^n$
for every color
$c\in C$,
because
$\frac{1}{p}<\frac{1}{p-1}+\frac{n+1-c}{n+1}<1-\frac{1}{p}$
due to the condition
$\frac{1}{p-1}+\frac{1}{p}<\frac{1}{n+1}$,
see sect.~\ref{subsect:targetree}. Therefore,
$\eta_0=\eta_0(c)+c\nu\md\Z^n$
lies in the interior of the cube
$c(A)=A+c\nu\md\Z^n$,
\end{proof}

It follows from Lemma~\ref{Lem:colorcov} that
$\dist(\eta_0,\d H^{-k}(c(A)))\to\infty$
as
$k\to -\infty$.
Therefore, given a vertex 
$v$
of
$Q_{c,k'}$, $k'\in\Z$,
the vertex 
$H^{-k}(c(A))$
of
$Q_{c,k}$
cover  
$v$
for all 
$k<0$
with sufficiently large
$|k|$.
The same argument as in the proof of Lemma~\ref{Lem:tree}
shows that
$T_c$
is a tree for every color
$c\in C$.

\begin{proof}[Proof of Proposition~\ref{Pro:covertree}] To complete
the proof of Proposition~\ref{Pro:covertree} it remains to show that
the set
$V_c$
of cubes in
$R^n$
satisfies the separation property for every color
$c\in C$.
That is for each
$k$, $k'\in\Z$, $k'<k$,
any vertex of
$V_{c,k}$
either is separated from every vertex of
$V_{c,k'}$
by the distance at least
$\la^{k+1}$
or it lies inside of some vertex of
$V_{c,k'}$
being separated from its boundary by the distance at least
$\la^{k+1}$, $\la=1/p$.

This property is already proved for the case
$k'\ge 0$.
Assume that
$k'<0$.
Then, by definition,
$V_{c,k'}=H^{-k'}(V_{c,0})$.
Note that
$H^{-1}=h_l:\R^n\to\R^n$
for the letter
$l=(1,\dots,1)\in L$.
Since
$h_w(V_{c,0})=V_{c,k}$
for each word
$w\in W_k$
by Lemma~\ref{Lem:shiftalph},
we have
$V_{c,k}=H^{-k}(V_{c,0})$
for all
$k\ge 0$
and, hence, for all
$k\in\Z$.
It follows
$$H^{-k'}(V_{c,k})=H^{k-k'}(V_{c,0})=V_{c,k-k'},$$
and the general case follows from the case
$k'\ge 0$.
\end{proof}

\subsection{Definition of the embedding
$f:\hyp^{n+1}\to\prod_{c\in C}T_c$}

It is convenient to rescale the metric of
$\hyp^{n+1}$
as follows. The space
$\R\times\R^n$
with the warped product metric
$$ds^2=dt^2+e^{2\si t}d\rho_n^2,$$
where
$d\rho_n^2=dx_1^2+\dots dx_n^2$
is the canonical Euclidean metric on
$\R^n$,
and
$\si=\ln p$,
has the constant curvature
$K\equiv -\si^2$.
In other words,
$\R\times\R^n$
with the metric
$ds^2$
is the hyperbolic space
$\hyp^{n+1}$
rescaled by the factor
$1/\si$,
$\hyp_p^{n+1}$
for brevity. Then the shift
$\pr:\R\times\R^n\to\R\times\R^n$,
$\pr(t,x)=(t-1,x)$,
is a homothety of
$ds^2$, 
restricted to any horosphere
${t}\times\R^n$,
with the coefficient
$\la=1/p$.

Now, for every color
$c\in C$
and every 
$k\in\Z$
we consider the copy
$(k,Q_{c,k})\sub\R\times\R^n$
of the set
$Q_{c,k}$,
see sect.~\ref{subsect:color}. Recall that
$Q_{c,k}\sub\R^n$
consists of cubes with diameter
$\la^ka\sqrt n$.
It follows that the diameter
of the cubes of
$(k,Q_{c,k})$
is
$a\sqrt n$
(w.r.t. the horospherical metric of
$(k,\R^n)$
induced by
$ds^2$),
i.e., it is one and the same for all levels
$k\in\Z$
and all colors
$c\in C$.

For every
$c\in C$,
we define a (discontinuous) map
$f_c:\hyp_p^{n+1}\to T_c$
assigning to 
$z\in\hyp_p^{n+1}$
the vertex
$f_c(z)\in V_c$
represented by a cube from
$\cup_k(k,Q_{c,k})$
closest to
$z$.
This defines the required map
$f:\hyp_p^{n+1}\to\prod_{c\in C}T_c$.

The following fact allows to simplify the proof of the Lipschitz 
property of
$f$.

\begin{Lem}\label{Lem:trisimple} Let
$a$, $b$, $c$
be the side lengths of a triangle in a metric space
such that
$c\ge a$.
Then
$a+b\le 3c$.
\end{Lem}

\begin{proof} If
$b<2a$
then
$a+b<3a\le 3c$.
Assume that
$b\ge 2a$.
Then
$2c\ge 2(b-a)\ge b$.
Therefore,
$3c\ge a+b$
as well.
\end{proof}

\begin{Pro}\label{Pro:lip} The map
$f_c:\hyp_p^{n+1}\to T_c$
is large scale Lipschitz for every
$c\in C$.
\end{Pro}

\begin{proof} We have to show that
$$\dist(f_c(z),f_c(z'))\le\La\dist(z,z')+\al$$
for some
$\La\ge 1$, $\al\ge 0$
and all
$z$, $z'\in\hyp_p^{n+1}$.
Fix
$z$, $z'\in\hyp_p^{n+1}$,
and put
$v=f_c(z)$, $v'=f_c(z')$. Then
$v\in V_{c,k}$, $v'\in V_{c,k'}$
for some
$k$, $k'\in\Z$.
W.L.G. we can assume that
$z$, $z'$
are the centers of the cubes
$v\sub(k,Q_{c,k})$, $v'\sub(k',Q_{c,k'})$.
respectively, and that
$k'\ge k$.
Taking the point
$z''\in (k',\R^n)$
which projects to 
$z\in (k,\R^n)$,
we note that
$\dist(f_c(z),f_c(z''))\le k'-k$
because the levels of the end point of any edge in
$T_c$
differ at least by 1, and
$k'-k=\dist(z,z'')$.
Using this and Lemma~\ref{Lem:trisimple} we can
assume W.L.G. that
$k'=k$,
i.e.,
the points
$z$, $z'$
belong to the horosphere
$(k,\R^n)\sub\hyp_p^{n+1}$,
and the cubes
$v$, $v'$
have one and the same level
$k$.

It follows from the definition of the edges that
no shortest path in
$T_c$
has an interior vertex with locally maximal level. Thus
the shortest path in
$T_c$
between
$v$
and
$v'$
has a unique vertex
$v_0$
of a lowest level
$k_0$, $k\ge k_0$.

Let
$v_1\in v_0v$, $v_1'\in v_0v'$
be the vertices adjacent to
$v_0$.
Note that
$v_1$, $v_1'\sub v_0$
considered as cubes in
$\R^n$.
Then by the separation property, the cubes
$v_1$, $v_1'$
either are disjoint or one of them is contained
in the other. However, the last case is excluded because
otherwise we would have a path in
$T_c$
between
$v_1$
and
$v_1'$
missing the vertex
$v_0$,
and hence the initial path
$vv_0\cup v_0v'$
would not be the shortest one.

Assuming that
$v_1\sub(k_1,Q_{c,k_1})$, $v_1'\sub(k_1',Q_{c,k_1'})$,
where
$k\ge k_1\ge k_1'>k_0$,
we obtain that the Euclidean distance between the cubes
$v_1$, $v_1'\sub\R^n$
is at least
$\la^{k_1+1}$.
Therefore, the horospherical distance between
$z$, $z'$
w.r.t. the horosphere
$(k,\R^n)$
is
$\ge e^{\si k}\la^{k_1+1}\ge\const(p)p^{k-k_1}$,
where
$\const(p)>0$
depends only on
$p$.
Then for the distance in
$\hyp_p^{n+1}$
we have
$$\dist(z,z')\ge\const_1(p)(k-k_1)-\const_2(p).$$

First, we consider the case that the next
vertex 
$v_2'$
of the path
$v_0v'$
following
$v_1'$
has the level
$k_2'\ge k_1$.
Then for the distances in the tree
$T_c$
we have
\[\dist(v_0,v)=1+\dist(v_1,v)\le 1+k-k_1\]
and
\[\dist(v_0,v')=2+\dist(v_2',v')\le 2+k-k_2'
  \le 2+k-k_1.\]
Therefore,
\[\dist(v,v')\le 3+2(k-k_1)\le\La\dist(z,z')+\al\]
for some constants
$\La>0$, $\al\ge 0$
depending only on
$p$.

Assume now that
$k_2'<k_1$.
Then we have 
$v'\sub v_2'\sub v_1'$
for the vertices
$v'$, $v_2'$, $v_1'$
considered as the cubes in
$\R^n$.
Thus the point
$z'$
projected down to the horosphere
$(k_1',\R^n)$
lies in the interior of the cube
$v_1'$
being separated from its boundary by the Euclidean distance
$\ge\la^{k_2'+1}$.
It follows that the horospherical distance between
$z$, $z'$ is
$\ge\la^{k_2'+1}e^{\si k}=p^{k-k_2'-1}$,
and consequently
$$\dist(z,z')\ge\const_1(p)(k-k_2')-\const_2(p).$$
On the other hand,
\[\dist(v_0,v)\le 1+k-k_1\le 1+k-k_2'\]
and
\[\dist(v_0,v')\le 2+k-k_2'.\]
Therefore,
\[\dist(v,v')\le 3+2(k-k_2')\le\La\dist(z,z')+\al\]
for some constants
$\La>0$, $\al\ge 0$
depending only on
$p$.
\end{proof}

The following Proposition completes the proof of Theorem~\ref{Thm:main1}.

\begin{Pro}\label{Pro:bilip} The map
$f:\hyp_p^{n+1}\to\prod_{c\in C}T_c$
is quasi-isometric.
\end{Pro}

\begin{proof} By Proposition~\ref{Pro:lip}, it remains to show that
$$\dist(z,z')\le\La\dist(f(z),f(z'))+\al$$
for some constants
$\La\ge 1$, $\al\ge 0$,
and all
$z$, $z'\in\hyp_p^{n+1}$.

We can assume that
$z=(k,x)$, $z'=(k',x')$
for some
$k$, $k'\in\Z$, $k'\ge k$,
where
$x$, $x'\in\R^n$.
First, consider the case
$x=x'$.
Then the geodesic segment
$zz'\in\hyp_p^{n+1}$
intersects
$k'-k+1$
horospheres
$(t,\R^n)$
at the points
$(k,x)$, $(k+1,x),\dots,(k',x)$.
Since
$\cup_{c\in C}V_{c,k}=\R^n$,
at least
$(k'-k+1)/|C|$
of these points belong to cubes with one
and the same color
$c\in C$.
All of those cubes contain the cube
$f_c(z')\in V_{c,k'}$.
Hence, for the distance in
$T_c$
we have
$\dist(f_c(z),f_c(z'))\ge (k'-k+1)/|C|-1$
by the separation property. Therefore,
$$\dist(f(z),f(z'))\ge\frac{1}{|C|}(k'-k+1)-1\ge\frac{1}{n+1}\dist(z,z')-1.$$

In general case, we consider the points
$\ov z$, $\ov z'$
which are projections of
$z$, $z'$
respectively to a horosphere
$(k_0,\R^n)$
with largest level
$k_0$,
for which the horospherical distance between
$\ov z$, $\ov z'$
is at most 
$p\sqrt n$. 
Then this distance is
$>\sqrt n$,
thus
$f_c(\ov z)\neq f_c(\ov z')$
for every color
$c\in C$.
Since the geodesics in every tree
$T_c$
have no interior point with locally maximal level,
it follows that
$$\dist(f_c(z),f_c(z'))\ge\dist(f_c(z),f_c(\ov z))+
   \dist(f_c(\ov z'),f_c(z')).$$
Applying the first case, we obtain
\begin{eqnarray*}
  \dist(f(z),f(z'))
  &\ge&\frac{1}{(n+1)}\max\{\dist(z,\ov z),\dist(z',\ov z')\}-1\\
  &\ge&\const_1(n)\dist(z,z')-\const_2(n,p)
\end{eqnarray*}
for some positive constants depending only on
$p$
and/or
$n$.
\end{proof}

\section{Proof of Theorem~\ref{Thm:main2}}

Actually, we prove a stronger result.

\begin{Thm}\label{Thm:nolargeball} Given
$l\ge 1$, $m\ge 0$, $n\in\N$,
there is
$r_0=r_0(l,m,n)>0$
such that no ball
$B_r\sub\hyp^2$
of radius
$r\ge r_0$
can be
$(l,m)$-quasi-isometrically
embedded in
$T\times\R^n$
for any tree
$T$.
\end{Thm}

Clearly, Theorem~\ref{Thm:main2} follows from 
Theorem~\ref{Thm:nolargeball}. In turn, Theorem~\ref{Thm:nolargeball}
is a corollary of the following more general result about
arcs with bounded winding. 

Let
$X$
be a CAT($-1$) space of bounded geometry. The last means
that there are
$\rho_X>0$
and
$M_X:(0,\infty)\to(0,\infty)$
such that every ball
$B_r\sub X$
of radius
$r>0$
contains at most
$M_X(r)$
points which are
$\rho_X$-separated.
We use notation
$|x-x'|$
for the distance in
$X$
between
$x$, $x'\in X$,
and
$\diam A$
for the diameter of
$A\sub X$
in
$X$.
We fix an origin 
$o\in X$,
and denote by
$S_r$
the metric sphere in
$X$
of radius
$r$
centered at
$o$.
If
$x$, $x'$
are different from
$o$,
we let
$\angle_o(x,x')$
be the angle at
$\ov o$
of the comparison triangle
$\ov o\ov x\ov x'\sub\hyp^2$.
If 
$A$
misses
$o$,
we put
$\angle_o(A)=\sup\set{\angle_o(x,x')}{$x,x'\in A$}$,
the angle diameter of
$A$.

Let
$A\sub S_r$, $r>0$,
be an arc. Any
$a$, $a'\in A$
define the subarc
$A(a,a')\sub A$
with the end points
$a$, $a'$.
Given
$\de$, $\si\in (0,1)$,
one says that the arc
$A$
has a
$\de$-bounded
winding at the scale
$\si$,
if for every subarc
$A(a,a')\sub A$
with
$\angle_o(A(a,a'))\ge\si\angle_o(A)$
we have
$|a-a'|\ge\de\diam A$
(this property is useful in the theory of quasi-conformal
mappings).

From now on, we assume that some constants
$l\ge 1$, $m\ge 0$
and an integer
$n\ge 0$
are fixed, and saying about a quasi-isometric map 
$f:A\to Y$
we mean
that the map 
$f$
is
$(l,m)$-quasi-isometric, i.e.,
$$\frac{1}{l}|a-a'|-m\le\dist(f(a),f(a'))\le l|a-a'|+m$$
for all
$a$, $a'\in A$.
The constants
$\si_0$, $r_0$, $N_0$,
which will be introduced below, depend, in particular,
on
$l$, $m$, $n$, $\rho_X$
and
$M_X$.
We do not reflect this dependence in notations for brevity.

\begin{Thm}\label{Thm:winding} For every
$\de\in(0,1/4]$, $\ep>0$
there are
$\si_0=\si_0(\de)\in(0,1)$, $r_0=r_0(\de,\ep)\ge 1$,
such that the following holds true. Let
$A\sub X$
be an arc with
\begin{itemize}
\item[(1)] $A\sub S_r$
for some
$r\ge r_0$;
\item[(2)] $\angle_o(A)\ge\ep$;
\item[(3)] $A$
has a 
$\de$-bounded
winding at the scale
$\si_0$.
\end{itemize}
Then there is no quasi-isometric map
$f:A\to T\times\R^n$
for any metric tree
$T$.
\end{Thm}

Any nondegenerate arc
$A\sub\d B_r\sub\hyp^2$
subtending the angle
$\le\pi$
has 1-bounded winding at any scale
$\si\in (0,1)$
(for
$\hyp^n$
with
$n\ge 3$
this is certainly not true). Hence, Theorem~\ref{Thm:nolargeball}
follows from Theorem~\ref{Thm:winding}.

\begin{Rem} The factor
$\R^n$
in Theorem~\ref{Thm:winding} can be replaced by any 
geodesic space
$Y$,
which satisfies the following condition. There
is a function
$N:(0,1)\to\N$
such that for every
$\rho\in (0,1)$
every ball 
$B_R\sub Y$
with sufficiently large radius
$R$
contains at most
$N(\rho)$
points which are
$\rho R$-separated
(the constants
$r_0$, $\si_0$
then depend also on
$N$).
\end{Rem}

\begin{Rem} Using Theorem~\ref{Thm:winding} one can show that
there is no quasi-isometric map
$f:X\to T\times\R^n$,
where
$X$
is a CAT($-1$)-space with bounded geometry such that
$\dim\di X=1$.
That is, the hyperbolic rank (see \cite{Gr}, \cite{BS1})
of the product
$T\times\R^n$
is zero for any tree
$T$
and any
$n\ge 0$.
\end{Rem}

\begin{Rem}\label{Rem:continuity} Every quasi-isometric
map 
$f:A\to T\times\R^n$
as above can be easily modified to a continuous one. So, 
W.L.G. we shall prove only that there is no {\em continuous} 
quasi-isometric map
$f:A\to T\times\R^n$.
\end{Rem}

Briefly, the proof proceeds as follows. Assuming that
the assertion is not true, 
we find, for sufficiently large
$r$,
an arc
$A_r\sub S_r$
with bounded winding at some scale, and a continuous 
quasi-isometric map
$f_r:A_r\to T\times\R^n$.
Let
$g_r:A_r\to T$
be the composition of
$f_r$
with the projection
$T\times\R^n\to T$
onto the first factor.
We study preimages
$g_r^{-1}(C)$
of geodesic segments
$C$
in the subtree
$D_r=g_r(A_r)\sub T$.
Since
$f_r$
maps
$g_r^{-1}(C)$
quasi-isometrically in
$C\times R^n\sub\R^{n+1}$,
the preimage
$g_r^{-1}(C)$
is small in the sense that it contains a bounded
amount of sufficiently separated points. Then it follows
from hyperbolicity of the space
$X$
that the complement
$D_r\sm C$
contains large subtrees, i.e., subtrees for 
which preimage has a large subarc in
$A_r$.
Since
$A_r$
has a bounded winding, the end points of such an arc
are sufficiently separated in 
$X$.
The key point of the proof (Lemma~\ref{Lem:manybig}) 
is that the number of large subtrees is sufficiently large 
for an appropriately chosen segment
$C\sub D_r$,
and hence there are sufficiently many separated points in
$A_r$
mapped by
$f_r$
in
$C\times\R^n$
to obtain a contradiction with properties of
$C\times\R^n$.

Assume that for
$r\ge 1$
there is an arc
$A_r\sub S_r$
and a continuous quasi-isometric map
$f_r:A_r\to T\times\R^n$.
Then
$D_r=g_r(A_r)\sub T$
is a connected compact subset and therefore it is a subtree.
 
\begin{Lem}\label{Lem:boundsep} Fix
$\de\in(0,1)$.
There is
$N_0=N_0(\de)\in\N$,
such that for every segment
$C\sub D_r$
the arc
$A_r$
contains at most
$N_0$
points from
$g_r^{-1}(C)$,
which are pairwise separated by the distance
$\ge\de\rho$
in
$X$,
where
$$\rho=\max\{\frac{\rho_X}{\de},\frac{2lm}{\de},\diam g_r^{-1}(C)\}.$$
\end{Lem}

\begin{proof} We put
$\de'=\frac{\de}{4l^2}$,
$r(l,m,\de)=2lm/\de$
and note that
$$r(l,m,\de)\ge m/l$$
since
$l\ge 1$
and
$\de<1$.
We define
$N_0$
as the maximum of two numbers
$N_0'$
and
$N_0''$,
where
$N_0'$
is the maximal number of
$\de'$-separated
points in the ball of radius 1 in
$\R^{n+1}$, $N_0'=N_0'(l,n,\de)$,
and 
$N_0''=\max\{M_X(\rho_X/\de),M_X(2lm/\de)\}$,
$N_0''=N_0''(l,m,X,\de)$.
Therefore,
$$N_0=N_0(l,m,n,X,\de)=N_0(\de)$$
according to our agreement.

If
$\diam g_r^{-1}(C)\le\rho$,
then the claim follows from the definition of
$N_0''$.
Thus we assume that
$\rho=\diam g_r^{-1}(C)>r(l,m,\de)$.
Consider a segment
$C\sub D_r$.
Let
$E\sub g_r^{-1}(C)$
be a maximal
$\de\rho$-separated
subset. Since
$\diam E\le\rho$,
the set
$f_r(E)$
lies in a ball of radius
$\le l\rho+m\le \rho'=2l\rho$
in
$C\times\R^n\sub\R^{n+1}$.
Furthermore,
$f_r(E)$
is
$(\frac{\de\rho}{l}-m)$-separated.
Since
$\frac{\de\rho}{l}-m\ge\frac{\de\rho}{2l}=\de'\rho'$,
we obtain
$|E|\le N_0$.
\end{proof}

Every segment
$C\sub D_r$
separates the tree
$D_r\sub T$
into a collection
$\cT(C)$
of closed subtrees in
$D_r$:
every subtree
$P\sub D_r$
from
$\cT(C)$
is the closure in
$D_r$
of some connected component of the complement
$D_r\sm C$.

Let
$P$
be a closed subtree in
$D_r$.
Every connected component of the preimage
$g_r^{-1}(P)\sub A_r$
is an arc, may be degenerate.
For
$\si\in(0,1)$,
the subtree
$P$
is called
$\si$-large
w.r.t.
$D_r$,
if
$g_r^{-1}(P)$
contains a connected component with the angle diameter
$\ge\si\cdot\angle_o(A_r)$
(every such a component is called large). If it is clear, which tree
$D_r$
is considered, then we speak about
$\si$-large
trees.

\begin{Lem}\label{Lem:existbig} Let
$\de\in(0,1/4]$, $N_0=N_0(\de)$
be the constant from Lemma~\ref{Lem:boundsep},
$\si\in(0,\frac{1}{2(N_0+1)}]$.
Then for every
$\ep>0$
there is
$r(\de,\ep)>0$,
such that for all
$r\ge r(\de,\ep)$
we have the following. If the arc
$A_r\sub S_r$
has the angle diameter
$\angle_o(A_r)\ge\ep$,
then for every segment
$C\sub D_r$,
the collection of trees
$\cT(C)$,
into which the segment separates the tree
$D_r$,
contains at least one
$\si$-large
tree.
\end{Lem}

\begin{proof} We fix
$\ep>0$
and take
$r(\de,\ep)>\max\{\frac{\rho_X}{\de},\frac{2lm}{\de}\}$
such that for all
$r\ge r(\de,\ep)$
the following holds true: if
$|x-x'|\le r$
for
$x$, $x'\in S_r$,
then the angle distance
$\angle_o(x,x')\le\frac{\ep}{2N_0}$.
Such an
$r(\de,\ep)$
exists since
$X$
is a CAT($-1$)-space.
Assume that
$r\ge r(\de,\ep)$
and
$\angle_o(A_r)\ge\ep$.
Then
$\rho=\diam A_r>\max\{\frac{\rho_X}{\de},\frac{2lm}{\de}\}$.

We choose a maximal
$\de\rho$-separated
subset
$E\sub g_r^{-1}(C)$.
The closed balls
$B_{\de\rho}(x)\sub X$, $x\in E$,
cover the preimage
$g_r^{-1}(C)$,
and the number of them
$|E|\le N_0$
according to Lemma~\ref{Lem:boundsep}, since
$$\diam g_r^{-1}(C)\le\rho.$$
Assuming that an orientation of the arc
$A_r$
is fixed, for every point
$x\in E$
we take the point
$x^+\in A_r$
of the first coming in the ball
$B_{\de\rho}(x)$
and the point
$x^-\in A_r$
of the last coming out from the ball
$B_{\de\rho}(x)$.
We have
$|x^--x^+|\le 2\de\rho\le r$
by the choice of
$\de\le1/4$
and because
$\rho\le 2r$.
Then the angle distance between these points satisfies
$\angle_o(x^-,x^+)\le\frac{\ep}{2N_0}$
by the choice of
$r$.

It is clear that
$$g_r^{-1}(C)\sub\bigcup_{x\in E}A(x^-,x^+).$$
The complement
$A_r\sm\cup_{x\in E}A(x^-,x^+)$
consists of open intervals, whose number
$\le|E|+1\le N_0+1$.
It suffices to prove that at least one of those intervals
has the angle diameter
$\ge\si\cdot\angle_o(A_r)$.
Assume that it is not the case, and every interval has
the angle diameter
$<\si\cdot\angle_o(A_r)$.
Then
$$\angle_o(A_r)<(N_0+1)\cdot\si\cdot\angle_o(A_r)+\sum_{x\in E}\angle_0(x^-,x^+)
  \le\frac{1}{2}(\angle_o(A_r)+\ep).$$
However, this contradicts the condition
$\angle_o(A_r)\ge\ep$.
\end{proof}

\begin{Lem}\label{Lem:manybig} Let
$\de\in(0,1/4]$,
and let
$N_0=N_0(\de)$
be the constant from Lemma~\ref{Lem:boundsep},
$\si\in(0,\frac{1}{2(N_0+1)}]$, $\ep>0$.
Then for every
$N\in\N$
and every
$r\ge r(\de,\si^N\ep)$
we have: if an arc
$A_r\sub S_r$
has the angle diameter
$\angle_o(A_r)\ge\ep$,
then there exists a segment
$C\sub D_r$
such that the tree collection
$\cT(C)$,
into which the segment separates the tree
$D_r$,
contains at least
$N$ $\si^N$-large
trees.
\end{Lem}

\begin{proof} The segment
$C$
will be constructed inductively as the union of an increasing
(in one direction) sequence of subsegments
$C_0\sub C_1\sub\dots$.
We take as
$C_0$
some extreme vertex of the tree
$D_r$.
Then the collection
$\cT(C_0)$
consists of one tree
$P_1=D_r$
which is, of course,
$\si$-large.
As
$C_1$
we take the edge of the tree
$P_1$
adjacent to
$C_0$.

Assume that we have already constructed segments
$C_0\sub C_1\sub\dots\sub C_{k-1}\sub D_r$
and
$\si$-large
trees
$P_1\in\cT(C_0),\dots,P_{k-1}\in\cT(C_{k-2})$
such that
$P_i$
is a unique
$\si$-large
tree in
$\cT(C_{i-1})$
for all
$i=1,\dots,k-1$,
and
$P_{k-1}\cap C_{k-1}$
is the edge of the tree
$P_{k-1}$
adjacent to the segment
$C_{k-2}$.
One can assume that
$$r(\de,\si^N\ep)\ge r(\de,\si^{N-1}\ep)\ge\dots
  \ge r(\de,\ep).$$
Then by Lemma~\ref{Lem:existbig}, the collection
$\cT(C_{k-1})$
contains at least one
$\si$-large
tree, and by the assumption of uniqueness, all such trees
of that collection have one and the same common point with the
segment
$C_{k-1}$,
which is its end (different from
$C_0$
if
$k\ge 2$).
If the collection
$\cT(C_{k-1})$
also contains a unique
$\si$-large
tree
$P_k$,
then we take
as
$C_k$
the union of the segment
$C_{k-1}$
and the edge of the tree
$P_k$
adjacent to that segment. By construction, the union
$C_k$
is a segment.

We assert that for some
$k\ge 1$
the collection
$\cT(C_k)$
contains at least two
$\si$-large
trees. Indeed, the tree
$D_r\sub T$
is a compact subset being the continuous image of the compact set
$A_r$, $D_r=g_r(A_r)$.
Thus it has only finitely many edges. If one assumes that
the assertion is not true, then the procedure described above gives
after a finite number of steps a segment
$C\sub D_r$,
that connects some extreme vertices of the tree
$D_r$,
such that the collection
$\cT(C)$
contains no
$\si$-large
trees. This contradicts Lemma~\ref{Lem:existbig}.

Let
$k\ge 1$
be the least integer for which the collection
$\cT(C_k)$
contains at least two
$\si$-large
trees. We denote by
$Q_1$
one of them, and for another,
$Q$,
we choose a large connected component
$A_1\sub A_r$
of its preimage
$g_r^{-1}(Q)$
and denote by
$D_1=g_r(A_1)$
the subtree
$D_1\sub Q$.
Then
$\angle_o(A_1)\ge\si\cdot\angle(A_r)\ge\si\ep$.
At least one of the ends of the arc
$A_1$
is an interior point in the arc
$A_r$
(otherwise
$A_1=A_r$
and
$D_1=D_r$,
which is impossible). Thus the tree
$D_1\sub Q$
has a common point with the segment
$C_k$
(which is an extreme point for both). The condition
$$r\ge r(\de,\si^N\ep)\ge\dots
  \ge r(\de,\si\ep)$$
allows to apply to
$D_1$
the previous arguments and continue the construction
of the segment
$C$,
starting with the end of the segment
$C_k$,
different from
$C_0$.
Every
$\si$-large
subtree
$P\sub D_1$
w.r.t.
$D_1$
is a
$\si^2$-large
w.r.t.
$D_r$
since for a large connected component
$A'$
of the set
$g_r^{-1}(P)\sub A_1\sub A_r$
its angle diameter
$\angle_o(A')\ge\si\cdot\angle_o(A_1)\ge\si^2\cdot\angle_o(A_r)$.

Since
$$r\ge r(\de,\si^N\ep),$$
the condition of Lemma~\ref{Lem:existbig}
is satisfied at least for
$N$
such steps, and we obtain a segment
$C\sub D_r$
and different trees
$Q_1,\dots,Q_N\in\cT(C)$,
where the tree
$Q_i$
is
$\si^i$-large
and hence
$\si^N$-large
w.r.t.
$D_r$, $i=1,\dots,N$.
\end{proof}

\begin{Lem}\label{Lem:boundbig} Let
$\de$, $\si\in(0,1)$, $r>0$.
Assume that an arc
$A_r\sub S_r$
has a
$\de$-bounded
winding at the scale
$\si$,
and
$\diam A_r\ge\max\{\frac{\rho_X}{\de},\frac{2lm}{\de}\}$.
Then for every segment
$C\sub D_r=g_r(A_r)$
the collection
$\cT(C)$
contains at most
$N_0+2$
subtrees which are
$\si$-large,
where
$N_0=N_0(\de)$
is the constant from Lemma~\ref{Lem:boundsep}.
\end{Lem}

\begin{proof} We fix a segment
$C\sub D_r$.
Let
$A'\sub A_r$
be a connected component of the set
$g_r^{-1}(P)$
for one of the trees
$P\in\cT(C)$.
Note that
$g_r(a)\in C$
for an end
$a$
of the arc
$A'$,
which is an interior point in
$A_r$.
This follows from continuity of the map
$g_r$.
It the tree
$P$
is large and
$A'$
is a large component from its preimage,
$\angle_o(A')\ge\si\cdot\angle_o(A_r)$,
then the ends
$a$, $a'$
of the arc
$A'$
are separated in
$X$
by the distance
$|a-a'|\ge\de\rho$,
where
$\rho=\diam A_r$,
according the property of a bounded winding of the arc
$A_r$.

Assume that trees
$P_1$, $P_2$, $P_3\sub D_r$
are large and there exist connected components
$A_1\sub g_r^{-1}(P_1)$,
$A_2\sub g_r^{-1}(P_2)$,
$A_3\sub g_r^{-1}(P_3)$,
such that the arc
$A_2$
separates the arcs
$A_1$
and
$A_3$
on the arc
$A_r$.
Then the distance in
$X$
between each end of the arc
$A_1$
and each end of the arc
$A_3$
is at least
$\de\rho$,
because the pairs of the points above are the ends of
subarcs in
$A_r$
containing
$A_2$
and hence having the angle diameter
$\ge\angle_o(A_2)\ge\si\cdot\angle_o(A_r)$.

Therefore, any maximal collection
$\cT_0$
of large trees from the collection
$\cT(C)$,
such that for every tree there is a large interior component
of the preimage in
$A_r$,
and each two of them are separated by a large component
of the preimage of some large tree (not necessarily from the
collection
$\cT_0$),
gives
$2k$
points from
$g_r^{-1}(C)$,
which are pairwise
$\de\rho$-separated
in
$X$,
where
$k$
is the number of the trees of the collection
$\cT_0$.
Since
$\rho\ge\max\{\frac{\rho_X}{\de},\frac{2lm}{\de},\diam g_r^{-1}(C)\}$,
we have
$2k\le N_0$
according Lemma~\ref{Lem:boundsep}.
On the other hand, it is clear that the number of the
large trees in
$\cT(C)$
is at most
$2k+2$
(the additional factor 2 takes into account large trees
which may not have large interior component in the preimage)
and thus that number is
$\le N_0+2$.
\end{proof}

\begin{proof}[Proof of Theorem~\ref{Thm:winding}] Let
$N_0=N_0(\de)$
be the constant from Lemma~\ref{Lem:boundsep}. We put
$$\si_0=(2(N_0+1))^{-(N_0+3)},\quad
  r_0=r(\de,\si_0\ep),$$
where
$r(\de,\ep)$
is the constant from Lemma~\ref{Lem:existbig}. Then
$\si_0=\si_0(\de)$
and
$r_0=r_0(\de,\ep)$.

Assume now that
$r\ge r_0$
and an arc
$A_r\sub S_r$
with the angle diameter
$\angle_o(A_r)\ge\ep$
is quasi-isometrically (and continuously) mapped into
$T\times\R^n$.
We let
$\si=\frac{1}{2(N_0+1)}$.
Then for
$N=N_0+3$
we have
$r\ge r(\de,\si^N\ep)$,
i.e., for the arc
$A_r$
the condition of Lemma~\ref{Lem:manybig} is satisfied, and
$\si^N=\si_0$.
According to that Lemma, there is a segment
$C\sub D_r=g_r(A_r)$
such that the collection of trees
$\cT(C)$
contains at least
$N$
trees
which are
$\si_0$-large.

Assume further that the arc
$A_r$
has a
$\de$-bounded
winding at the scale
$\si_0$.
By Lemma~\ref{Lem:boundbig}, the collection
$\cT(C)$
contains at most
$N_0+2$
trees
which are
$\si_0$-large.
This is a contradiction with what we get above, since
$N_0+2<N$.
Therefore, the arc
$A_r$
cannot be quasi-isometrically mapped into
$T\times\R^n$.
\end{proof}


\bigskip
\begin{tabbing}

Sergei Buyalo,\hskip11em\relax \= Viktor Schroeder,\\ 

St. Petersburg Dept. of Steklov \>
Institut f\"ur Mathematik, Universit\"at \\

Math. Institute RAS, Fontanka 27, \>
Z\"urich, Winterthurer Strasse 190, \\

191023 St. Petersburg, Russia\>  CH-8057 Z\"urich, Switzerland\\

{\tt buyalo@pdmi.ras.ru}\> {\tt vschroed@math.unizh.ch}\\

\end{tabbing}

\end{document}